\documentclass[12pt,a4paper]{article}
\usepackage[utf8]{inputenc}
\usepackage{amsmath}
\usepackage{amsfonts}
\usepackage{amssymb}
\usepackage{amsthm}
\usepackage{authblk}
\usepackage{inputenc}
\begin{document}
\title{A better method than t-free for Robin's hypothesis}
\author{Xiaolong Wu}
\affil{Ex. Institute of Mathematics, Chinese Academy of Sciences}
\affil{xwu622@comcast.net}
\date{November 29, 2018}
\maketitle
\begin{abstract}
Robin made hypothesis that $\sigma (n)<e^\gamma n \log \log n$ for all integers $n>5040$. For a positive integer $t>1$, an integer N is called t-free if N cannot be divided by prime power $p^t$ for any prime p. Some works have shown if N is t-free, then N satisfies Robin's inequality, for t=5, 7, 11, 16. This article shows that the condition of t-free can be reduced to "N cannot be divided by $2^t$". I proved that if N cannot be divided by $2^{17}$, then N satisfies Robin’s inequality.
\end{abstract}

Robin made a hypothesis [Robin 1984] that the Robin's inequality 
\begin{equation}\tag{RI}
\sigma (n)<e^\gamma n \log \log n, 
\end{equation}
 holds for all integers $n>5040$. Here $\sigma (n)=\sum_{d|n} d$ is the divisor sum function, $\gamma $ is the Euler-Mascheroni constant, log is the nature logarithm.\\
\indent Robin also proved that his hypothesis is true if and only if Riemann hypothesis is true.

 For calculation convenience, we define \[\rho (n):=\frac{\sigma (n)}{n}.\]
Then Robin’s inequality can also be written as 
\begin{equation}\tag{RI}
\rho (n)<e^\gamma \log \log n .
\end{equation}

    Let $N>5040$ be an integer. Write the factorization of N as \[N=\prod_{i=1}^{r} p_i^{a_i},\]
where r is the number of distinct prime factors of N, $p_i$ are listed in increasing orders. So $p_r$ is the largest prime factor of N.

     According to $[Morrill; Platt\, 2018]$, (RI) holds for all integers  $n, 5040<n\leq 10^{(10^{13})}$. So, we assume $N>10^{(10^{13})}$.

    For a positive integer $t>1$, an integer N is called t-free if N cannot be divided by prime power $p^t$ for any prime p. Some works have shown if N is t-free, then N satisfies Robin's inequality, for t=5 {[}CLMS 2007{]}, 7 {[}Solé; Planat 2012{]}, 11 {[}Broughan;Trudgian 2015{]}, 16 $[Morrill; Platt\, 2018]$.  $[Morrill; Platt\, 2018]$ also shows 25-free integers satisfy (RI). But that is by a different method of numerical calculation as that article said "We calculated that $R_{25}(N_{177\, 244\, 758\, 016})<1-10^{-16}$".
    
    According to abundant number theory, any possible counter-example of (RI) has its exponents of primes in decreasing order, so, the condition "t-free" is basically equivalent to exponent of prime 2 is less than t.
    
    The previous works on t-free used $\zeta(t)$. However, we have $2^{-t}<\zeta(t)<2^{-t+1}$ for integer $t\geq 2$, and \[\lim_{t\rightarrow \infty}2^t\zeta(t)=1,\quad  t\,runs\,over\,integers\geq 2.\] Hence, using $\zeta(t)$ has almost the same effect as using $2^{-t}$.
    
    This article shows that the condition of t-free can be reduced to "N cannot be divided by $2^t$". In particular,  I proved: if N cannot be divided by $2^{17}$, then N satisfies Robin’s inequality.

    Using a method similar to abundant number theory, I can improve the result to $2^{42}$. I will post that improvement later.\\[8mm]
\noindent {\bfseries Lemma 1. (Mertens’ third theorem) }
\textit{ For any integer $n>7\,713\,133\,853$, we have
\begin{equation}\tag{L2.1}
\sum_{p\leq n}\log \left(\frac{p}{p-1}\right)=\log \log n+\gamma+R(n),
\end{equation}
where $\gamma$ is the Euler-Mascheroni constant, R(n) is the remainder such that
\begin{equation}\tag{L2.2}
 -\frac{0.005586}{(\log n)^2}<R(n)<\frac{0.005586}{(\log n)^2}.
\end{equation}
}
\begin{proof}
Just follow the proof of Theorem 5.6 of [Dusart 2018] with $k=2,\eta_2=0.01$.
\begin{align*}
\lvert R(n)\lvert &<\frac{0.01}{2(\log n)^2}+\frac{4}{3}\cdot \frac{0.01}{(\log n)^3}\\
&=\frac{0.01}{(\log n)^2}\left(\frac{1}{2}+\frac{4}{3\log n}\right)<\frac{0.005586}{(\log n)^2}.\tag{L2.3}
\end{align*}
\end{proof}

\textbf{Definition.} Let N be an integer, p be its largest prime factor. We define a \textbf{generalized factorization} of N as the product
\[N=\prod_{p_i\leq p}p_i^{a_i}\] 
where $p_i$ runs through all primes $p_i\leq p$. If a prime $p_i<p$ does not divide N, we set $a_i=0$, i.e. $p_i^{a_i}=p_i^0=1$ if $p_i$ is not a factor of N.\\[5mm]
\noindent {\bfseries Theorem 1. }
\textit{
Let $N>10^{(10^{13})}$ be an integer. Define
\[T:=\frac{0.005586}{\left(\log \log \left(10^{(10^{13})}\right)\right)^2}, \quad  E(p):=\left\lfloor \frac{-\log T}{\log p}\right\rfloor,\]
here $\lfloor x \rfloor$ denotes the largest integer $\leq x$. Write the generalized  factorization of N as
\[N=\prod_{p_i\leq p_r}p_i^{a_i},\]
where $p_i$ are listed in increasing order and $p_r$ is the largest prime factor of N.
If $p_s^{E(p_s)}$ does not divide N for some $p_s$, then N satisfies (RI).\\
\indent In particular, since
\[\log T=-12.04,  \quad E(2)=\left\lfloor \frac{12.04}{\log 2}\right\rfloor=17,\]
if $2^{17}$ does not divide N, then N satisfies (RI).}
\begin{proof}
We use induction by assuming that all integers n in $(5040, N)$ with exponent of $p_s$ less than $E(p_s)$ satisfy (RI). We know
\begin{equation}\tag{1.1}
\rho (N)=\prod_{p_i\leq p_r}\frac{p_i-p_i^{-a_i}}{p_i-1}=\prod_{p_i\leq p_r}\frac{1-p_i^{-a_i-1}}{1-p_i^{-1}}.
\end{equation}
Define 
\begin{equation}\tag{1.2}
P(N):=\prod_{p_i\leq p_r}\frac{p_i}{p_i-1}=\prod_{p_i\leq p_r}\frac{1}{1-p_i^{-1}}.
\end{equation}
Define
\begin{align*}
S(N):&=\log P(N)-\log \rho (N)\\
&=-\sum_{p_i\leq p_r}\log (1-p_i^{-a_i-1})=\sum_{p_i\leq p_r}\sum_{k=1}^{\infty}\frac{1}{kp_i^{(a_i+1)k}},\tag{1.3}
\end{align*}
Then apply Lemma 1 to $\log P(N)$,
\begin{equation}\tag{1.4}
\log \rho (N)=\log P(N)-S(N)<\log \log p_r+\gamma+R(p_r)-S(N).
\end{equation}
where
\begin{equation}\tag{1.5}
R(p_r)=\frac{0.005586}{(\log p_r)^2}.
\end{equation}
By method of Lemma 11 of {[}CLMS 2007{]} or Theorem 26 of {[}NY 2014{]}, we may assume $p_r<\log N$. 

Since $\log \log x+\gamma+R(x)$ is an increasing function in x, we may substitute $p_r$ with $\log N$ in (1.4),
\begin{equation}\tag{1.6}
\log \rho (N)<\log \log \log N+\gamma+R(\log N)-S(N),
\end{equation}

    Now, 
\begin{equation}\tag{1.7}
S(N)>\sum_{p_i\leq p_r}p_s^{-a_s-1}\geq p_s^{-a_s-1}\geq p_s^{-E(p_s)}.
\end{equation}
Combine (1.5) and (1.7), we have
\begin{align*}
R(\log N)-S(N)&<\frac{0.005586}{(\log \log (10^{(10^{13})}))^2}-p_s^{-E(p_s )}\\
&\leq T-p_s^{\frac{\log T}{\log p_s}}=T-e^{(\log p_s)\left(\frac{\log T}{\log p_s}\right)}=T-T=0.\tag{1.8}
\end{align*}
Then from (1.6) we get $\log \rho (N)<\log \log \log N+\gamma$. Take exponential and we are done.
\end{proof}

\noindent {\bfseries Corollary 1. }
\textit{
Let $N>10^{(10^{13})}$ be an integer. Define
\[L:=\prod_{p\leq p_r}p^{E(p)}.\]
If N is not divisible by L, then N satisfies (RI).}
\begin{proof}
If N is not divisible by L, there must exist a prime power, say $p_s^{E(p_s)}$, of L that cannot divide N. Then N satisfies (RI) by Theorem 1.                                                 
\end{proof}

{\bfseries Remark:} The size of L can be estimated as $L\approx 10^{(10^{4.868})}$. 
\begin{center}
{\bfseries \large References}
\end{center}

\noindent {[}Broughan;Trudgian 2015{]} K.A. Broughan and T. Trudgian. \textit{Robin’s inequality for 11-free integers}. Integers, 15: Paper No. A12, 5, 2015.\\
{[}CLMS 2007{]} Y.-J. Choie, N. Lichiardopol, P. Moree, and P. Sol\'{e}. \textit{On Robin’s criterion for the Riemann hypothesis}. J. Th\'{e}or. Nombres Bordeaux, 19(2):357–372, 2007.\\
{[}Dusart 2018{]} P. Dusart. \textit{Explicit estimates of some functions over primes}. Ramanujan J.,
45(1):227–251, 2018.\\
{[}Morrill;Platt 2018{]} T. Morrill, D. Platt. \textit{Robin’s inequality for 25-free integers and obstacles to analytic improvement} \\
https://arxiv.org/abs/1809.10813\\
{[}NY 2014{]} S. Nazardonyavi and S. Yakubovich. \textit{Extremely Abundant Numbers and the Riemann Hypothesis} Journal of Integer Sequences, Vol. 17 (2014),Article 14.2.8.\\
{[}Robin 1984{]} G. Robin. \textit{Grandes valeurs de la fonction somme des diviseurs et hypoth\'{e}se de Riemann}. Journal de mathématiques pures et appliqu\'{e}es. (9), 63(2):187–213, 1984.\\
{[}Solé; Planat 2012{]} P. Sol\'{e} and M. Planat. \textit{The Robin inequality for 7-free integers}. Integers, 12(2):301–309, 2012.

\end{document}